\documentstyle[11pt]{article}

\newfont{\Bb}{msbm10 scaled\magstephalf}

\begin{document}
 \noindent

 \begin{center}
 {\LARGE Periodic motion of a charge  on a manifold in the  magnetic fields}
 \end{center}
 \vspace{10mm}
 \begin{center}
  \noindent
    {\large  Guangcun Lu}\footnote{Partially supported by the NNSF
   19501021  of China.}\\[5pt]
       Nankai Institute of Mathematics, Nankai University\\
               Tianjin 300071, P. R. China\\[5pt]  \vspace{-2mm}
               (E-mail: lugc@public.tjuc.com.cn)\\[4mm]

               Preliminary version\hspace{2mm}May 18, 1999\\[5pt]
              \vspace{-4mm}
              Revised version\hspace{2mm} May 20, 2000  \\[5pt]
              \end{center}

  \vspace{10mm}
 \begin{abstract}
\baselineskip 15pt
 In this paper we prove the existence of a periodic motion of a charge
 on  a large class of manifolds under the action of the 
 magnetic fields. Our methods also give a class of closed manifolds whose
cotangent bundles contain no the closed exact Lagrangian submanifolds.

 \end{abstract}

MSC 1991: 58F05 58E05

\baselineskip=18pt

\section{ Introduction and results}

\subsection{The question}

The periodic motion question of a charge on an Riemannian manifold $(N, g)$ in
the magnetic field( abbreviated to ``PMMQ" below) is a very important and 
difficult question in the mathematics and physics([Ar1][No]). It can be formulated as
 
    {\bf PMMQ}.\quad Looking for the nonconstant periodic solutions of Hamiltonian 
system  $$\dot z=X_{H_g}(z)$$  on the energy level $E_c=\{H_g=c\}$ with $c>0$, where 
$H_g:T^\ast N\to \mbox{\Bb R}$ is given by $H_g(z)=\frac{1}{2}\|z\|^2_g$
and $X_{H_g}$ is the Hamiltionian vector field of $H_g$ with respect to
the twisted symplectic form $\omega=\omega_{\rm can}+\pi_N^\ast\Omega$
on $T^\ast N$, the closed $2$-form $\Omega$ on $N$ corresponds to the magnetic
field.
 
In order to study it S.P.Novikov invented the variational principle of multi-valued
functionals([No][GN][NT][T]), V.I.Arnold introduced the symplectic topology methods(
[Ar2][Gi1]). On the detalied arguments of the history and progress
of this question  before 1995 the readers may refer to Ginzburg's beautiful survey 
paper [Gi1]. In addition, as showed by Example 3.7 in [Gi1] or Example 4.2 in [Gi2]
one cannot expect that the above question has always a solution. Thus it becomes 
very important to study some conditions under which the above question holds.

\subsection{The symplectic topology methods}

A well-known question in symplectic geometry is Weinstein conjecture, which claims:
every hypersurface $S$ of contact type in a symplectic manifold $(Q,\omega)$ carries
at least a closed characteristic([W]). Here $S$ is said to be of contact type if
there exists a transversal vectorfield $X$ defined on some open neighbrhood $U$ of
$S$ such that $L_X\omega=\omega$([W]). It is not difficult to check that
the energy level $E_c=\{H_g=c\}$ above cannot be of contact type in the sense if the
magnetic field $\Omega$ is not exact( see Remark 1.6.A in \S1.6 below).
Thus PMMQ is different from the Weinstein conjecture in the symplectic
manifold $(T^\ast N,\omega_{\rm can}+\pi_N^\ast\Omega)$. Fortunately, motivated
by the study of the latter Hofer and Zehnder introduced an important notion,
 Hofer-Zehnder symplectic capacity, which can not only be used to study the
Weinstein conjecture but also PMMQ above. Let us first recall it.
 Given a symplectic manifold
$(Q,\omega)$ we denote by ${\cal H}(Q,\omega)$ the subset of
$C^\infty(Q,\mbox{\Bb R})$ consisting of all smooth functionas with the following
properties:
\begin{itemize}
\item There exists a compact subset $K\subset Q\setminus\partial Q$ depending on $H$
such that $H|_{Q\setminus K}\equiv m(H)$, i.e. a constant;
\item There exists a nonempty open subset $U$ depending on $H$ such that
$H|_U\equiv 0$;
\item $0\le H(x)\le m(H)$ for all $x\in Q$.
\end{itemize}
We shall call $H\in{\cal H}(Q,\omega)$ admissible( resp. C-admissible) if it
has the property that all $T$-periodic( resp. contractible) solutions of
$\dot x=X_H(x)$ on $Q$ having periods $0<T\le 1$ are constant solutions.
Writting ${\cal H}_{ad}(Q,\omega)$( resp. ${\cal H}_{cad}(Q,\omega)$) the
sets of admissible ( resp. C-admissible) $H\in{\cal H}(Q,\omega)$, we define
$$C_{HZ}(Q,\omega)=\sup\{ m(H)\,|\,H\in{\cal H}_{ad}(Q,\omega)\},\quad
\bar C_{HZ}(Q,\omega)=\sup\{ m(H)\,|\,H\in{\cal H}_{cad}(Q,\omega)\}.$$
Obviously, both are symplectic invariants, and it always holds that
$C_{HZ}(Q,\omega)\le \bar C_{HZ}(Q,\omega)$. $\bar C_{HZ}$ was first introduced
in [Lu1]. We still call it Hofer-Zehnder symplectic capacity.
As proved by Hofer and Zehnder([HZ])
and strengthened by Struwe ([St]) if $C_{HZ}(Q,\omega)<+\infty$ then for any
compact hypersurface $S\subset Q$ and any embedding $\imath :S\times [0, 1]\to Q$
there exists a set $J\subset [0, 1]$ of parameters of measure $1$ such that
for every $s\in J$ the Hamiltonian flow on $\imath(S\times\{s\})$ carries a periodic
orbit. We have showed in [Lu1][Lu2] that the similar conclusion for $\bar C_{HZ}$ 
still holds. More precisely saying, if $\bar C_{HZ}(\imath(S\times [0, 1]),\omega)$ 
is finite then for every parameter $s$ in a set $J$ of measure $1$ the Hamiltonian
 flow on $\imath(S\times\{s\})$ carries a contractible (in $\imath(S\times\{s\})$ )
 periodic orbit since $\imath(S\times\{s\})$ has the same homotopy type as
$\imath(S\times [0, 1])$. However there exists a difference between $C_{HZ}$ and
$\bar C_{HZ}$. That is, it only satisfies the following monotonicity axiom of 
weaker form:
  
\noindent{\it For a symplectic embedding $\psi:(M_1,\omega_1)\to (M_2,\omega_2)$
  of codimension zero, if either $M_1$ is simply connected or $\psi$
  induces an injective homomorphism $\psi_\ast:\pi_1(M_1)\to\pi_2(M_2)$
  then it holds that
  $${\bar C}_{HZ}(M_1,\omega_1)\le {\bar C}_{HZ}(M_2,\omega_2).$$}\vspace{-2mm}
  
Hence if we can prove that for a given $c>0$ there exists a sufficiently small number
$\epsilon>0$ such that $C_{HZ}(\{c-\epsilon\le H_g\le c+\epsilon\},\omega)<+\infty$(
resp. $\bar C_{HZ}(\{c-\epsilon\le H_g\le c+\epsilon\},\omega)<+\infty$) then for
generic $c^\prime$ near $c$ the level $E_{c^\prime}$ carries a (resp.
contractible in $E_{c^\prime}$ ) nonconstant periodic orbit  of $X_{H_g}$.

\subsection{Some recent results}

Before 1995, for the case of the nonexact magnetic field(i.e. $\Omega$ a nonexact form)
on the high dimension the general result is little.
For the rational magnetic field, that is, a closed $2$-form $\Omega$ on $N$ satisfying
$$m(N,\Omega)={\inf}\{\langle[\Omega],\alpha\rangle>0\,|\,\alpha\in\pi_2(N)\}>0,
\eqno(1.1)$$
which is called the rationality index of $\Omega$, we obtained a general result
in Corollary E.2 of [Lu1]. Precisely speaking, for any rational closed two-form 
$\Omega$ on a closed smooth manifold $M$ and any Riemannian metric $g$ on 
$N:=M\times\mbox{\Bb R}/{2\pi\mbox{\Bb Z}}$ 
we proved that 
$$C_{HZ}(\{v\in T^\ast N\,|\,\|v\|_g\le c\}, 
\omega_{\rm can}+\pi_N^\ast(p_1^\ast\Omega))<+\infty$$
for sufficiently small $c>0$, where $p_1:N\to M$ is the natural
projection. This implies that there exists a nonconstant periodic orbit
of Hamiltonian vector field $X_{H_g}$ on the energy level $E_c$
for almost all (resp. sufficiently small) $c>0$( in the sense of measure theory)
if $m(M, \Omega)=+\infty$(resp. $0<m(M,\Omega)<+\infty$). 
For the periodic motions on the tours under the action of magnetic fields
$\Omega$ Mei-Yue Jiang proved that the generic levels of $H_g$ carry a 
nonconstant periodic orbit provided that de Rham cohomology class $[\Omega]$
is rational ([J1]). Recently, Ginzburg and Kerman removed the rationality
assumption on $[\Omega]$ in ([GiK]). It is very surprising that  L. Polterovich 
used Hofer's geometry approach to show that  if for any nonzero magnetic field 
$\Omega$ and  Riemannian metric $g$ on $\mbox{\Bb T}^n$  there exists a sequence
of positive energy values $c_k=c_k(g,\Omega)\to 0$ such that every level $\{H_g=c_k\}$
carries a nonconstant {\it contractible} closed orbit([P2]). For the case 
of the exact  magnetic field readers may refer to [Vi2] [BT].

\subsection{The exact Lagrangian embedding and normal submanifolds}
 
 Recall that a submanifold $L$ of middle dimension in a 
 symplectic manifold $(Q,\omega)$ is called {\it normal} if there is
 a field of Lagrangian subspaces along $L$ which is transversal to $L$([Si]). 
 In [P1] it was proved that each Lagrangian submanifold $L$, and those submanifolds
 which are sufficiently $C^1$-close to $L$, and each parallelizable
 totally real submanifold of $Q$ with respect to some 
 $J\in{\cal J}(Q,\omega)$ are all normal. 

\noindent{\bf Example 1.4.A.}\hspace{2mm}Every closed orientable $3$-dimensional 
manifold  is parallelizable totally real submanifold of $\mbox{\Bb C}^3$, and thus 
a normal  submanifold in $(\mbox{\Bb C}^3,\omega_0)$. It should be noted that
 $\mbox{\Bb S}^3$ is such a manifold and satisfies: 
$H^1(\mbox{\Bb S}^3,\mbox{\Bb R})= H^2(\mbox{\Bb S}^3,\mbox{\Bb R})=0$,
 but it can not be embedded into $(\mbox{\Bb C}^3, \omega_0)$
 in the Lagrangian way because there is no any closed simply connected Lagrangian
 submanifold in $(\mbox{\Bb C}^n, \omega_0)$.  On the other hand, 
 the necessary and sufficient condition of $n$-dimensional totally
 real closed submanifold in $\mbox{\Bb C}^n$ was obtained in [A]( also refer to
 [Th 3.2.4, ALP]).

The following proposition is a key to proof of Theorem 1.5.A.
We believe itself to have some independent importance.

\noindent{\bf Proposition 1.4.B.}\hspace{2mm}{\it Let $L$ be a closed normal
submanifold in a symplectic manifold $(Q,\sigma)$, and $g$
a Riemannian metric on $L$. Assume that $\sigma$ is exact near $L$, i.e.,
$\sigma=d\tau$ for some one-form near $L$.
Then for any $c>0$ there exists a $0<\delta_0(c)<1$
such that for every $0<\delta\le\delta_0(c)$ we have a symplectic embedding of
codimension zero $E_\delta$ from $(\{H_g\le c\}, \omega^L_{\rm can}+
\pi_L^\ast(\sigma|_L))$ into $(Q,\frac{1}{\delta}\sigma)$. Consequently,
for some $0<\delta_1<\delta_0$ and all $0<\delta\le\delta_1$ there also exists
a symplectic embedding of codimension zero from $(\{H_g\le c\}, \omega^L_{\rm can})$
into $(Q,\frac{1}{\delta}\sigma)$.  Moreover, for a
given open neighborhood of $L$ in $Q$ we can require $\delta>0$ so small that
the images of these symplectic embeddings are contained in this open neighborhood.}

From Gromov's striking theorem that {\it there is no exact Lagrangian embedding of
a closed manifold into $(\mbox{\Bb R}^{2l},\omega_0)$}([Gr]) we immediately obtain
the following corollary.

\noindent{\bf Corollary 1.4.C.}\hspace{2mm}{\it For any closed normal submanifold 
in a symplectic manifold $(\mbox{\Bb R}^{2l},\omega_0)$ there  is no exact Lagrangian
embedding of a closed manifold into $T^\ast L$.}

For more results on the generalization of Gromov theorem the readers may refer to
[V2][V3].

The monotonicity of the symplectic capacity directly leads to
 
\noindent{\bf Corollary 1.4.D.}\hspace{2mm}{\it Let $L$ be a closed normal
submanifold in a symplectic manifold $(Q,\sigma)$. If there exists 
a neighborhood of $L$, ${\cal U}$ such that $\sigma$ is exact on it and
$({\cal U},\sigma)$ has finite symplectic capacity. Then for any Riemannian metric
$g$ on $L$ and every $c>0$ the symplectic manifold $(\{H_g\le c\}, \omega^L_{can})$
has finite capacity.}

\subsection{Main results}

The manifolds in our main results below will be the product manifold $N=M\times L$ 
of a closed smooth manifold $M$ and a compact normal submanifold $L$ without
boundary of  $(\mbox{\Bb R}^{2l}, \omega_0)$. Denote by $P_M:N\to M$ is the natural
projection to the second factor and by 
$P_M^\ast: H_{de}^2(M;\mbox{\Bb R})\to H_{de}^2(N;\mbox{\Bb R})$ the  
homomorphism between their second de Rham cohomology groups induced by it. 
We shall omit the subscript ``de" in the de Rham cohomology groups below.
Then $P^\ast_M(H^2(M,\mbox{\Bb R}))$ is a subspace of $H^2(N,\mbox{\Bb R}))$.
Similarily, for any differomorphism $\phi\in{\rm Diff}(N)$ we denote by
$\phi^\ast$ the induced isomorphism between $H^2(N;\mbox{\Bb R})$ and itself.
We get a subset of $H^2(N;\mbox{\Bb R})$ as follows
$$\bigcup_{\phi\in{\rm Diff}(N)}\phi^\ast(P^\ast_M(H^2(M,\mbox{\Bb R}))).\eqno(1.2)$$
It seems to be strange. However, if $H^2(L,\mbox{\Bb R})=0$ and 
$H^1(L,\mbox{\Bb R})=0$ it directly follows from
the K\"unneth formula that $H^2(N,\mbox{\Bb R})$ and $H^2(M,\mbox{\Bb R})$ is 
isomorphic. Thus in the case the set in (1.2) is equals to $H^2(N,\mbox{\Bb R})$.

 Our first result to PMMQ above is

\noindent{\bf Theorem 1.5.A.}\hspace{2mm}{\it Let $N=M\times L$ be 
a closed smooth manifold $M$ and a compact normal submanifold $L$ without boundary 
of  $(\mbox{\Bb R}^{2l}, \omega_0)$ as above. $\Omega$ is a  closed two-form
on $N$  whose de Rham cohomology class $[\Omega]$ belongs to the set in (1.2). If
$\Omega|_{\pi_2(N)}=0$, then for every Riemannian metric $g$ on $N$ and $c>0$ it
holds that
$$C_{HZ}(\{H_g\le c\}, \omega)<+\infty,$$
where $\omega:=\omega^N_{\rm can}+\pi^\ast_N\Omega$. Consequently, for
generic $c>0$ the level $E_c=\{H_g=c\}$ carries a nonconstant periodic  orbit of
$X_{H_g}$. Here $X_{H_g}$ is the Hamiltonian vector field determined
by $i_{X_{H_g}}\omega=dH_g$. Specially, if $L$ is simply connected we can also
guarantee such  generic levels  $E_c=\{H_g=c\}$ to carry a nonconstant periodic 
orbit with the contractible projection to $N$.}

\noindent{\bf Corollary 1.5.B.}\hspace{2mm}{\it If $H^2(L,\mbox{\Bb R})=0$ and 
$\pi_1(M)$ is a  finite group then for any Riemannian metric $g$ and a 
closed $2$-form $\Omega$ on $N$ with $\Omega|_{\pi_2(N)}=0$ the generic level
$E_c=\{H_g=c\}$ carries a nonconstant periodic  orbit of
$X_{H_g}$.}

If a  submanifold $L$ of $(\mbox{\Bb R}^{2l},\omega_0)$ is Lagrangian, rather
than normal Theorem 1.5.A can be strengthened. The following is the second result 
to PMMQ.

\noindent{\bf Theorem 1.5.C.}\hspace{2mm}{\it Let $N=M\times L$ be a product
of a closed smooth manifold $M$ and a closed Lagrangian submanifold $L$ of 
$(\mbox{\Bb R}^{2l}, \omega_0)$.  $\Omega$ is a rational closed two-form on $N$  
whose de Rham cohomology class $[\Omega]$ belongs to the set in (1.2).
Then for every Riemannian metric $g$ on $N$ 
 there is an upper semi-continuous function 
 $\Gamma^g_\Omega: [0, +\infty)\to [0, +\infty)$ ( see (3.13))
such that for every $c>0$ with $\Gamma^g_\Omega(c)< \sqrt{m(N,
\Omega)/\pi}$ it holds that
$$\liminf_{\epsilon\to 0+}C_{HZ}(U(g,c,\epsilon), \omega)<\pi(\Gamma^g_\Omega(c))^2<
 m(N,\Omega),\eqno(1.3)$$
where $U(g, c,\epsilon):=\{z\in T^\ast N\,|\, c-\epsilon\le H_g(c)\le c+\epsilon\}$.
Consequently, for almost all $c^\prime>0$ near $c$ the level
$E_{c^\prime}=\{H_g=c^\prime\}$ carries a nonconstant periodic  orbit of 
$X_{H_g}$, where $X_{H_g}$ is the Hamiltonian vector
field of $H_g$ with respect to the symplectic form
$\omega=\omega_{\rm can}+\pi^\ast_N\Omega$.}
 
 \noindent{\bf Remark 1.5.D.}\hspace{2mm}In Corollary 1.5.C, if we denote by 
 $$c(\Omega, g):=\inf_{(\phi,\widehat\Omega, \alpha)}\|\alpha\|_g,$$
 where $\|\alpha\|_g=\sup_{z\in N}\sqrt{g(\alpha(z),\alpha(z))}$ and the
 infimum is taken over all possible  $(\phi,\widehat\Omega, \alpha)$ satisfying (2.8).
 Then when $c(\Omega, g)$ is small enough and $m(N,\Omega)$ is large enough the 
inequality
 $$\Gamma^g_\Omega(c)< \sqrt{m(N,\Omega)/\pi}\eqno(1.4)$$
 always holds for sufficiently small $c>0$.

 In fact, for any given sufficiently small $\varepsilon>0$ we may choose
 $(\phi,\widehat\Omega, \alpha)$
 such that $\|\alpha\|_g<c(\Omega, g)+\varepsilon$ and the image set of 
$\alpha$( as a section) is contained
 in the image of $\Upsilon$ in (3.4). This is possible if $c(\Omega,g)$ is small
 enough. Notice that  in this case $\Theta$ is a symplectomorphism from
 $(T^\ast N,\omega)$ to the symplectic manifold in (3.3)
 with $P^\ast_L(\lambda_0|_L)=0$.
  Hence for sufficiently small $c>0$ the set
 $\Theta(\{H_g=c\})$ is also contained in the image of $\Upsilon$. Denote by
 $$r(c, g,\Theta,\Upsilon):=\inf\{r>0\,|\, P_{\cal U}[(\Upsilon|_{{\rm Im}(\Upsilon)})^{-1}\circ\Theta(\{H_g=c\})]
\subseteq{\cal U}\cap F(Z^{2l}(r)),\,F\in{\rm Symp}(\mbox{\Bb R}^{2l},\omega_0)\},$$
 then when $m(N,\Omega)$ is so large that 
$r(c, g,\Theta,\Upsilon)<\sqrt{m(N,\Omega)/\pi}$, or more general
${\cal U}\subset Z^{2l}(r_0)$ for some $0<r_0<\sqrt{m(N,\Omega)/\pi}$ it holds that
 $$\Gamma^g_\Omega(c)<\sqrt{m(N,\Omega)/\pi}.$$ 
   
In some special cases we can get better results. For example, let
$L=\mbox{\Bb S}^1=\mbox{\Bb R}/2\pi\mbox{\Bb Z}$ and $\Omega$ 
 a rational closed $2$-form on $N=M\times\mbox{\Bb S}^1$ whose de Rham cohomology
class belongs to the set in (1.2).  For any Riemannian metric $g$ on $N$
we define the function $\Xi^g_\Omega: [0, +\infty)\to [0, +\infty)$  by
  $$\Xi^g_{\Omega}(c):={\inf}_{(\Theta,\widehat\Omega,\epsilon)}
\Bigl|\int_{\bar P_S(\Theta((U(g,c,\epsilon)))}
\omega^S_{\rm can}\Bigr|,\eqno(1.5)$$
where the infimum is taken over all pairs $(\Theta,\widehat\Omega)$
satisfying Lemma 2.3 and all $\epsilon>0$.

 Our third result to PMMQ is given as follows. 

 \noindent{\bf Theorem 1.5.E.}\hspace{2mm}{\it Under the above assumptions, if
$c>0$ is such that $\Xi^g_\Omega(c)< m(N,\Omega)$ then 
$$\liminf_{\epsilon\to 0+}C_{HZ}(U(g,c,\epsilon), \omega)<\pi(\Xi^g_\Omega(c))^2<
 m(N,\Omega),$$
and therefore for almost all $c^\prime>0$ near $c$
the levels $\{H_g=c^\prime\}$ carries a nonconstant periodic  orbit of $X_{H_g}$, 
where $X_{H_g}$ is the Hamiltonian vector field of $X_{H_g}$ with respect to the 
symplectic form $\omega=\omega_{\rm can}+\pi^\ast_N\Omega$.
In addition, the function $\Xi^g_\Omega$ is  upper semi-continuous, and
 also satisfies:
$$\Xi^g_{\lambda\Omega}(\lambda c)=\lambda \Xi^g_\Omega(c)\quad
\forall c>0,\;\lambda>0.\eqno(1.6)$$}

 This result can be generalized to the case that $L=T^n$. 

\noindent{\bf Corollary 1.5.F.}\hspace{2mm}{\it Let $\Omega$ be
 a rational closed $2$-form on $N=M\times T^n$ whose de Rham cohomology
class belongs to the set in (1.2). Then for any Riemannian metric $g$ on $N$ 
there exists a nonnegative upper semi-continuous
function $\widehat\Xi^g_\Omega: [0, +\infty)\to [0, +\infty)$ such that
for every $c>0$ with $\widehat\Xi^g_\Omega(c)<m(N,\Omega)$
and almost all $c^\prime>0$ near $c$
the levels $\{H_g=c^\prime\}$ carries a nonconstant periodic  orbit of $X_{H_g}$, 
where $X_{H_g}$ is the Hamiltonian vector field of $H_g$ with respect to the 
symplectic form $\omega=\omega_{\rm can}+\pi^\ast_N\Omega$.
 Specially, the function
$\widehat\Xi^g_\Omega$ also satisfies
$$\widehat\Xi^g_{\lambda\Omega}(\lambda c)=\lambda \widehat\Xi^g_\Omega(c)
\quad\forall c>0,\;\lambda>0.$$}

With the method of the proof of Corollary 1.5.B and Theorems 1.5.C, 1.5.E we can
easily arrive at the following corollary.

 \noindent{\bf Corollary 1.5.G.}\hspace{2mm}{\it In Theorem 1.5.C, 1.5.E,
if $H^2(L,\mbox{\Bb R})=0$ and $\pi_1(M)$ is a  finite group then for any Riemannian 
metric $g$ and any rational closed $2$-form $\Omega$ on $N$ the corresponding
conclusions therein hold.}

\noindent{\bf Remark 1.5.H.}\hspace{2mm}
One may think that it is difficult to determine the values of
  $\Gamma^g_\Omega$ and $\Xi^g_\Omega$ in the  Theorems and Corollaries above.
   But we affirm them to be finite numbers, and if $\Omega|_{\pi_2(N)}=0$
 then $m(N,\Omega)=+\infty$ and thus the conditions of these theorems are
 always satisfied in this case for all $c>0$. On the other hand 
since the functions  $\Gamma^g_\Omega$ and $\Xi^g_\Omega$ are upper semi-continous 
the sets  $\{c>0\,|\,\Gamma^g_\Omega(c)<\sqrt{m(N,\Omega)/\pi}\}$
 and $\{c>0\,|\,\Xi^g_\Omega(c)<m(N,\Omega)\}$ are open.
 Hence if they are not empty then there exist nonconstant periodic
 solutions on the levels $\{H_g=c\}$ for all $c$ in a positive measure set([St]).
On the another hand one may think that it is difficult to understand the meaning
of functions $\Gamma^g_\Omega$ and $\Xi^g_\Omega$. Our starting points are to attempt
using a series of symplectic embeddings of codimension zero to reduce our
question to the case for which Theorem 2.1 may be applied, and to guarantee each step
being optimal so that Theorem 2.1 is best applied. Both functions are to 
characterize the optimization in the way of the quantity. If the rationality
condition $m(M,\omega)>0$ in Theorem 2.1 can be removed then our arguments
show that the generic energy levels carry a nonconstant
Hamiltonian periodic orbit. However, it is regrettable for us not to be able
to remove this assumption yet.

Our final result to PMMQ is about the case of tours $T^n$.
 Let $\Omega$ be a magnetic 
field (a closed $2$-form) on it, and $\omega=\omega_{\rm can}+\pi_{T^n}^\ast\Omega$.
As pointed out in \S1.3 one had known that for any metric $g$ on $T^n$ generic
levels $E_c=\{H_g=c\}$ in $(T^\ast T^n, \omega)$ carries a nonconstant periodic
orbit of $X_{H_g}$([GiK]). When $\Omega$ is not exact we can furthermore obtain

\noindent{\bf Theorem 1.5.I.}\hspace{2mm}{\it If $n\ge 2$ and the de Rham cohomology
 class $[\Omega]$ is nonzero then for generic $c>0$ the levels $E_c=\{H_g=c\}$
 carries a nonconstant periodic orbit of $X_{H_g}$ whose projection to the base
$T^n$ is {\rm contractible}.}

\subsection{Two remarks}

\noindent{\bf Remark 1.6.A.}\hspace{2mm}For a given Riemannian metric $g$ on $N$ 
and $c>0$ we denote by
$$\Sigma^g_c:=\{v\in T^\ast N\,|\, \|v\|_g=c\}.$$ 
Let a closed two-form $\Omega$ on $T^\ast N$ be such that
$\omega_{\rm can}+\Omega$ is a symplectic form on $T^\ast N$.
One may ask whether for sufficiently large $c>0$ the hypersuface $\Sigma_c^g$
is of contact type in the symplectic manifolds
$(T^\ast N,\omega_{\rm can}+\Omega)$? If this holds then our partial results
may be derived from one in [Vi2]. The following proposition answers this
question.

\noindent{\bf Proposition 1.6.}\hspace{3mm}{\it If $\Omega$ is not exact 
the hypersuface $\Sigma_c^g$ cannot be of contact type in the symplectic manifolds
$(T^\ast N,\omega_{\rm can}+\Omega)$.}

\noindent{\it Proof}.\hspace{3mm}Assume $\Sigma^g_c$ to be of contact type for
$\omega_{\rm can}+\Omega$, then there exists a sufficiently small $\epsilon>0$ 
such that $\omega_{\rm can}+\Omega$ is exact on 
$U(g,c,\epsilon):=\{v\in T^\ast N\,|\,c-\epsilon<\|v\|_g<c+\epsilon\}$.
Thus there is a one-form $\alpha$ on $U(g,c,\epsilon)$ such that
$$\Omega=d\alpha\;\;{\rm on}\;\;U(g, c,\epsilon)$$
because $\omega_{\rm can}$ is always exact.
Notice that $\pi_N:T^\ast N\to N$ induces a natural  isomorphism 
$\pi_N^\ast:H^\ast(N,\mbox{\Bb R})\to H^\ast(T^\ast N,\mbox{\Bb R})$   
and $\pi_N^\ast([\Omega|_N])=[\pi_N^\ast(\Omega|_N)]=[\Omega]$, we have
$\pi_N^\ast(\Omega|_N)-\Omega=d\beta$ for some one-form $\beta$ on $T^\ast N$.
Hence 
$$\pi_N^\ast(\Omega|_N)=d(\alpha+\beta)\;\;{\rm on}\;\;U(g, c,\epsilon)$$
Take a smooth section $F:N\to T^\ast N$ of $\pi_N$ such that
$F(N)\subset U(g,c,\epsilon)$(e.g., $F$ is the zero section) then as a smooth map from $N$ to $T^\ast N$ it
satisfies:
$$F^\ast\circ\pi_N^\ast(\Omega|_N)=F^\ast d(\alpha+\beta)=d(F^\ast(\alpha+\beta)).$$
But $F^\ast\circ\pi_N^\ast(\Omega|_N)=(\pi_N\circ F)^\ast(\Omega|_N)$ and
$\pi_N\circ F=id_N$. Hence we get
$$\Omega|_N=d(F^\ast(\alpha+\beta)).$$
That is, $\Omega|_N$ is a exact form on $N$. Denote by $\gamma:=F^\ast(\alpha+\beta)$.
Then 
$$\Omega=\pi_N^\ast(\Omega|_N)-d\beta=\pi_N^\ast(d\gamma)-d\beta=
d(\pi_N^\ast\gamma-\beta),$$
this leads to a contradiction.\hfill$\Box$\vspace{2mm}

\noindent{\bf Remark 1.6.B.}\hspace{2mm}Our results always deal with
manifolds of product forms . These are due to the  limitation of our methods.
The following example shows that increasing a factor manifold in the base
manifold will have a real influence on periodic orbits of magnetic fields.

Let $M$ be a compact surface equipped with a metric $g_0$ of constant curvature
$K=-1$ and $\Omega$ the area form on $M$. From Example 3.7 in [Gi1] we know
that if $c>1$ on the level $E_c$ there are no closed characteristic with
contractible projections with respect to the symplectic structure
$\omega_M=d\lambda_M +\pi_M^\ast\Omega$.  Consider $N:=M\times\mbox{\Bb S}^3$
and $\omega=d\lambda_N +\pi_N^\ast(p_M^\ast\Omega)$. Here $p_M:N\to M$
is the natural projection. Notice that $m(M,\Omega)=+\infty$. By Corollary 1.6.C,
for any Riemannian metric $g_1$ and the product metric $g=g_0\times g_1$ on $N$
and generic $c>0$ the levels $E_c:=\{H_g=c\}$ carries a nontrivial closed
characteristic with the contractible projection to $N$. On the other hand,
Example 3.7 in [Gi1] showed that for $c>1$ the Hamiltonian flow of $X_{H_{g_0}}$
with respect to $\omega_M$ on $E^M_c:=\{H_{g_0}=c\}$ has no any closed 
characteristic with contractible projections to $M$.  Notice that
$E^M_c\times 0(\mbox{\Bb S}^3)\subset E_c$ for the zero section $0(\mbox{\Bb S}^3)$
of $T^\ast\mbox{\Bb S}^3$. 

  Our arguments are the symplectic topology methods.    
 In \S2 we give some lemmas and prove Proposition 1.4.B. The proofs of all theorems 
and  corollaries are given in \S3. Finally, some concluding remarks are given in \S4.

{\bf Acknowledgements}\hspace{2mm}
This work was done under Professor Claude Viterbo's guidance while the author was 
visiting IHES. I would like to express my hearty thanks to him for his many valuable 
suggestions and often sending his preprints to me. I wish to acknowledge 
Professor J. Bourguignon for his invitation and hospitality. I also thank Professor
L. Polterovich for his kind comments and detalied improving suggestions on 
the contents and writting style of this paper, and Professors V.I. Ginzburg 
and A. Taimanov for sending me many papers.

 \section{Some lemmas and proof of Proposition 1.4.B} 
 
Let us first recall  the following theorem.

\noindent{\bf Theorem 2.1}([Th C, Lu1]).\hspace{2mm}{\it Let $(M,\omega)$ 
be a strong geometrical bounded symplectic manifold with 
$m(M,\omega)>0$. If $r\in (0, \sqrt{m(M,\omega)/\pi})$ then
$${\bar C}_{HZ}(M\times Z^{2n}(r), \omega\oplus\omega_0)\le\pi r^2,
\eqno(2.1)$$
where $Z^{2n}(r):=\{(x_1,y_1,\cdots,x_n,y_n)\in\mbox{\Bb R}^{2n}\,|\,
x^2_1 + y^2_1\le r^2\}$.}

For $C_{HZ}$ this result was obtained in [FHV][HV] for $M$ closed, and
in [Ma] for $(M,\omega)=(T^\ast Q,\omega_{\rm can})$.

In view of our result mentioned in \S1.3, a natural question is what 
conclusions one can get for a general
closed two-form $\Omega$ on $N$ rather than on $M$. It was  Professor
Claude Viterbo who asked whether for a given closed two-form on $N$
 there exist a closed two-form $\widehat\Omega$ on
$M$ and a $\Psi\in{\rm Diff}(T^\ast N)$ such that
$$\omega_{\rm can}+\pi_N^\ast\Omega=\Psi^\ast(\omega_{\rm can}+
\pi_N^\ast(P_M^\ast\widehat\Omega + P_L^\ast(\omega_0|_L)))\,?\eqno(2.2)$$
This idea motivates the studies of this paper.

The following lemma directly follows from the local coordinate arguments. 

\noindent{\bf Lemma 2.2.}\quad{\it For any closed $2$-form 
$\Omega$ and $1$-form $\alpha$ on a manifold $N$ the diffeomorphism 
$\Psi:T^\ast N\to T^\ast N$ given by
$$(m,v)\mapsto (m, v+ \alpha(m))\eqno(2.3)$$
satisfies: 
$$\Psi^\ast(\omega_{\rm can}+\pi_N^\ast\Omega)=
\omega_{\rm can}+\pi_N^\ast\Omega+\pi^\ast_N(d\alpha).\eqno(2.4)$$
That is, $\Psi$ is a symplectomorphism from $(T^\ast N, \omega_1)$
to $(T^\ast N, \omega)$, where
$$\omega:=\omega_{\rm can}+\pi_N^\ast\Omega\quad{\rm and}\quad
\omega_1:=\omega_{\rm can}+\pi_N^\ast\Omega+\pi^\ast_N(d\alpha).\eqno(2.5)$$}

\noindent{\bf Lemma 2.3}\hspace{2mm}{\it Let $N=M\times L$ be as in Theorem 1.5.A.
If the de Rham cohomology class $[\Omega]$ of a given closed two-form $\Omega$ on 
$N$  belongs to the set in (1.2) then there exists a closed two-form 
$\widehat\Omega$ on $M$ with
$$m(N,\Omega)=m(M,\widehat\Omega)\eqno(2.6)$$
such that $(T^\ast N,\omega)$ is symplectomorphic to $(T^\ast N,\widehat\omega)$,
where $\widehat\omega$ is given by
$$\widehat\omega:=\omega_{\rm can}+\pi_N^\ast(P_M^\ast\widehat\Omega 
+ P_L^\ast(\omega_0|_L)).\,\eqno(2.7)$$
Specially, if $H^2(L, \mbox{\Bb R})=0$ and either $H^1(M, \mbox{\Bb R})=0$ or 
$H^1(L,\mbox{\Bb R})=0$
then for every closed two-form $\Omega$ on $N$ the above conclusions hold.}

\noindent{\it Proof.}\hspace{2mm}Since $[\Omega]$  belongs to the set in (1.2)
there exists a $\phi\in{\rm Diff}(N)$ such that
$[\phi^\ast\Omega]=\phi^\ast[\Omega]$ belongs to
$P_M^\ast(H^2(M,\mbox{\Bb R}))$, and thus there must exist a closed two-form 
$\widehat\Omega$ on $M$ and a one-form $\alpha$ on $N$ such that 
$$\phi^\ast\Omega-P_M^\ast\widehat\Omega=d\alpha.\eqno(2.8)$$
But such $\phi$, $\widehat\Omega$ and $\alpha$ may not be
unique.  Notice that $\phi$ may lift to a symplectomorphism 
$\Phi:(T^\ast N,\omega_{\rm can})\to
(T^\ast N, \omega_{\rm can})$ by the formula
$$\Phi(m,v)=(\phi(m), [d\phi(m)^{-1}]^\ast(v)).\eqno(2.9)$$
They satisfy that $\pi_N\circ\Phi=\phi\circ\pi_N: T^\ast N\to N$
and therfore
$$\Phi^\ast\circ\pi_N^\ast=\pi_N^\ast\circ\phi^\ast:\Omega^2(N)\to
\Omega^2(T^\ast N).$$
By the definition of $\omega$ in (2.5) we get
$$\Phi^\ast\omega=\omega_{\rm can}+\pi^\ast_N(\phi^\ast\Omega)=
\omega_{\rm can}+ \pi_N^\ast(P_M^\ast\widehat\Omega)+
d(\pi_N^\ast\alpha)=\widehat\omega +d(\pi_N^\ast (\alpha-
P_L^\ast(\lambda_0|_L)))\eqno(2.10)$$
because $\omega_0=d\lambda_0$ is the standard symplectic form
on $\mbox{\Bb R}^{2l}$.
By Lemma 2.2 there exists a diffeomorphism 
$\Psi\in{\rm Diff}(T^\ast N)$ given by
$$(m,v)\mapsto (m, v+ \alpha(m)-P_L^\ast(\lambda_0|_L)(m))
\eqno(2.11)$$
such that $\Psi^\ast\widehat\omega=\Phi^\ast\omega$ or
$(\Psi\circ\Phi^{-1})^\ast\widehat\omega=\omega$.
Then
$$\Theta=\Theta(\Omega,\widehat\Omega, \phi,\alpha):=\Psi\circ\Phi^{-1}\eqno(2.12)$$
 given by
$$(m, v)\mapsto\Bigl(\phi^{-1}(m), d\phi(m)^\ast(v)+ 
(\alpha-P_L^\ast(\lambda_0|_L))(\phi^{-1}(m))
\Bigr)\eqno(2.13)$$
is a symplectomorphicism from $(T^\ast N,\omega)$ to $(T^\ast N,\widehat\omega)$.
As to (2.6), notice that $P_M$ induces a surjective homomorphism 
$P_{M\ast}:\pi_2(N)\to\pi_2(M)$
and that $\phi$ induces an isomorphism $\phi_\ast:\pi_2(N)\to\pi_2(N)$,
it may follow from (1.1) and the equalities
\begin{eqnarray*}
\langle[\Omega],\beta\rangle=\langle (\phi^{-1})^\ast\circ\phi^\ast([\Omega]),
      \beta\rangle
&=&\langle \phi^\ast([\Omega]),\phi_\ast^{-1}(\beta)\rangle\\
&=&\langle P^\ast_M([\widehat\Omega]), (\phi^{-1})_\ast(\beta)\rangle=
\langle[\widehat\Omega], P_{M\ast}(\phi^{-1}_\ast(\beta))\rangle\;\,
\forall\beta\in\pi_2(N).
\end{eqnarray*}
The final claim is a direct
consequence of K\"unneth formula because in this case 
$P_M$ induces an isomorphism $P_M^\ast:H^2(M, \mbox{\Bb R})\to H^2(N, \mbox{\Bb R})$.
Consequently, Lemma 2.3 holds.\hfill$\Box$\vspace{2mm}

In order to prove Theorem 1.5.A we also need the following lemma, which
perhaps goes back to the early work of Weinstein and Givental.

\noindent{\bf Lemma 2.4}(cf.[P1, Prop.1.9] and [Si, Th. 2.4]).\hspace{2mm}{\it
For a closed normal submanifold $L$ in a symplectic manifold $(Q,\sigma)$
and the restriction $\sigma|_L$ there exists an open  neighbourhood
${\cal U}$ of $L$ in $Q$ and an embedding $\varphi:{\cal U}\to T^\ast L$
such that $\varphi(L)$ coincides with the zero section of $T^\ast L$
and $\varphi^\ast(\omega^L_{\rm can}+ \pi^\ast_L(\sigma|_L))=\sigma$.}

For a closed smooth manifold $L$ of dimension $l$ we choose an atals
$\{(U_\alpha,\alpha)\}$ consisting of $m$ local coordinate charts. We
also require each $\alpha(U_\alpha)$ to be equal to the unit ball
$B^l$ centred at origin in $\mbox{\Bb R}^l$. For $q\in U_\alpha$ let
$\alpha(q)=(x^\alpha_1(q),\cdots, x^\alpha_l(q))$. It induces an obvious
bundle trivialization $\Phi_\alpha: T^\ast U_\alpha\to B^l\times\mbox{\Bb R}^l$
given by
$$(q, v^\ast)\mapsto (x^\alpha_1(q),\cdots, x^\alpha_l(q); y^\alpha_1(q, v^\ast),
\cdots, y^\alpha_l(q, v^\ast)),$$
where $y^\alpha_i(q, v^\ast)$ are determined by
$v^\ast\circ d\pi_L(q, v^\ast)=\sum^l_{j=1}y^\alpha_j(q, v^\ast) dx^\alpha_j(q)$.
The following lemma is very key to the proof of Theorem 1.5.A.

\noindent{\bf Lemma 2.5.}\hspace{2mm}{\it Let $g$ be a Riemannian metric on $L$
and $\lambda$ a $1$-form on $L$. Then for given  positive numbers $a$ and
$\epsilon\le 1$ there exist $b_\epsilon>0$ and a smooth function
$K_{a,\epsilon}: T^\ast L\to [0,1]$ such that
\begin{description}
\item[(i)] $K_{a,\epsilon}(v^\ast)=1$ for $\|v^\ast\|_g\le a$, and
$K_{a,\epsilon}(v^\ast)=0$ for $\|v^\ast\|_g\ge b_\epsilon$;
\item[(ii)]  for $K^\alpha_{a,\epsilon}:=(\Phi_\alpha^{-1})^\ast K_{a,\epsilon}:
B^l\times\mbox{\Bb R}^l\to [0,1]$ and $\lambda^\alpha=(\alpha^{-1})^\ast\lambda=
\sum^l_{i=1}\lambda_i^\alpha dx_i^\alpha$ it holds that
$$\biggl|\sum^l_{i=1}\frac{\partial K^\alpha_{a,\epsilon}}{\partial y^\alpha_i}
\lambda^\alpha_i\biggr|<\epsilon;$$
\item[(iii)] $A\omega_{\rm can}^L + A d(\pi_L^\ast\lambda) + 
t(1-A)d(K_{a,\epsilon}\pi_L^\ast\lambda)$ are symplectic
forms on $T^\ast L$ for every $A>1$ and $t\in [-1,1]$.
\end{description}}

\noindent{\it Proof.}\hspace{2mm}Choose a smooth function $\gamma: [0,\infty)\to
[0,1]$ such that (i) $\gamma(t)=1$ for $t\le 1$, (ii) $\gamma(t)=0$ for $t\ge 3$,
(iii) $|\gamma^\prime(t)|\le 1$. For each $\varepsilon>0$ we define
$$\gamma_{a,\varepsilon}(t):=\gamma(\varepsilon t+ 1-a\varepsilon).$$
Then $\gamma_{a,\varepsilon}(t)=1$ for $t\le a$,
$\gamma_{a,\varepsilon}(t)=0$ for $t\ge a+ 2/\varepsilon$,
 and $|\gamma_{a,\varepsilon}^\prime(t)|\le\varepsilon$ for all $\varepsilon>0$.
Denote by $H_{a,\varepsilon}((q, v^\ast))=\gamma_{a,\varepsilon}(\|v^\ast\|_g)$.
It is clearly smooth. Moreover, the local expression of it, 
$H_{a,\varepsilon}^\alpha:=(\Phi_\alpha^{-1})^\ast H_{a,\varepsilon}$, is given by
$$H_{a,\varepsilon}^\alpha(x^\alpha_1,\cdots, x^\alpha_l; y^\alpha_1,\cdots,
 y^\alpha_l)=\gamma_{a,\varepsilon}\biggl(\sqrt{\sum^l_{i, j}g^\alpha_{ij}(x)
y^\alpha_iy^\alpha_j}\biggr).$$ 
Thus
$$\frac{\partial H_{a,\varepsilon}^\alpha}{\partial y^\alpha_i}(x^\alpha, y^\alpha)=
\gamma_{a,\varepsilon}^\prime\biggl(\sqrt{\sum^l_{i, j}g^\alpha_{ij}(x^\alpha)
y^\alpha_i y^\alpha_j}\biggr)\frac{\sum^l_{j=1}g^\alpha_{ij}(x^\alpha)y^\alpha_j}
{\sqrt{\sum^l_{i, j}g^\alpha_{ij}(x^\alpha)y^\alpha_i y^\alpha_j}}.\eqno(2.14)$$
Denote by
$$c(\lambda):=\max\bigl\{|\lambda^\alpha_i(x)|\,\bigm|\,x\in\alpha(U_\alpha),
\,1\le i\le l,\,1\le\alpha\le m\bigr\}.$$
Notice that there exists a constant $C(g)>0$ such that
$$\frac{\sum^l_{i,j}|g^\alpha_{ij}(x^\alpha)||y^\alpha_j|}{\sqrt{\sum^l_{i, j}
g^\alpha_{ij}(x^\alpha)y^\alpha_i y^\alpha_j}}\le C(g)$$
for all $x^\alpha\in\alpha(U_\alpha)$ and $y^\alpha$ and $1\le\alpha\le m$.
We get
$$\Bigl|\sum^l_{i=1}\frac{\partial H_{a,\varepsilon}^\alpha}{\partial y^\alpha_i}
(x^\alpha, y^\alpha)\lambda_i^\alpha(x^\alpha)\Bigr|\le\sum^l_{i=1}
\Bigl|\frac{\partial H_{a,\varepsilon}^\alpha}{\partial y^\alpha_i}
(x^\alpha, y^\alpha)\Bigr||\lambda^\alpha_i(x^\alpha)|\le \varepsilon c(\lambda) C(g)
\eqno(2.15)$$
for all $x^\alpha\in\alpha(U_\alpha)$ and $y^\alpha$ and $1\le\alpha\le m$.
Hence it suffices to choose $\varepsilon_0=\varepsilon_0(\epsilon)>0$ such that 
$$\varepsilon_0 C(g) c(\lambda)<\epsilon.$$ 
Then $K_{a,\epsilon}:=H_{a,\varepsilon_0}$ satisfies (i)(ii).  

As to (iii), note that $A\omega_{\rm can}^L + A d(\pi_L^\ast\lambda) + 
t(1-A)d(K_{a,\epsilon}\pi_L^\ast\lambda)$ has the
following local expression
\begin{eqnarray*}
&&A \Bigl[\sum^l_{i=1}dx^\alpha_i\wedge dy^\alpha_i
+ \sum^l_{i,j}\frac{\partial\lambda^\alpha_j}{\partial x^\alpha_i}dx^\alpha_i\wedge
dx^\alpha_j\Bigr]+\\
&&t(1-A)\Bigl[ H^\alpha_{a,\varepsilon_0}
\sum^l_{i,j}\frac{\partial\lambda^\alpha_j}{\partial x^\alpha_i}dx^\alpha_i
\wedge dx^\alpha_j
+\sum^l_{i,j}\frac{\partial H^\alpha_{a,\varepsilon_0}}{\partial x^\alpha_i}
\lambda^\alpha_j dx^\alpha_i\wedge dx^\alpha_j-
\sum^l_{i,j}\frac{\partial H^\alpha_{a,\varepsilon_0}}{\partial y^\alpha_j}
\lambda^\alpha_i dx^\alpha_i\wedge dy^\alpha_i\Bigr],
\end{eqnarray*}
whose matrix in the natural basis $\partial/\partial x_1^\alpha,\cdots
 \partial/\partial x_l^\alpha;\partial/\partial y_1^\alpha,\cdots
 \partial/\partial y_l^\alpha$ is 
$S=\left(\begin{array}{cc}
S_{xx}&S_{xy}\\
 -S^t_{xy}&S_{yy}
\end{array}\right)$. Here $S_{yy}=0$, $S_{xx}=(a_{ij})$ and 
$S_{xy}=AI_l-t(1-A)(b_{ij})$. The matrix elements $a_{ij}$ and $b_{ij}$ are given by
$$a_{ij}=A\biggl(\frac{\partial\lambda^\alpha_j}{\partial x^\alpha_i}-
\frac{\partial\lambda^\alpha_i}{\partial x^\alpha_j}\biggr)+ t(1-A)\biggl[
H^\alpha_{a,\varepsilon_0}(\frac{\partial\lambda^\alpha_j}{\partial x^\alpha_i}-
\frac{\partial\lambda^\alpha_i}{\partial x^\alpha_j})+
 \frac{\partial H^\alpha_{a,\varepsilon_0}}{\partial x^\alpha_i}
\lambda^\alpha_j- \frac{\partial H^\alpha_{a,\varepsilon_0}}{\partial x^\alpha_j}
\lambda^\alpha_i\biggr],\eqno(2.16)$$\vspace{-2mm}
$$b_{ij}=\frac{\partial H^\alpha_{a,\varepsilon_0}}
{\partial y^\alpha_j}\lambda^\alpha_i.\eqno(2.17)$$
Notice that $S$ is nonsingular if and only if $S_{xy}$ is so. 
Assume that $S_{xy}\zeta=0$ for some vector $\zeta=(\zeta_1,\cdots,\zeta_l)^t$
in $\mbox{\Bb R}^{l}$. Then it holds that
$$A\zeta=t(1-A)(\lambda^\alpha_1,\cdots, \lambda^\alpha_l)^t \Bigl(
\frac{\partial H^\alpha_{a,\varepsilon_0}}{\partial y^\alpha_1},\cdots,
\frac{\partial H^\alpha_{a,\varepsilon_0}}{\partial y^\alpha_l}\Bigr)\zeta=
t(1-A)\Bigl(\sum^l_{i=1}\frac{\partial H^\alpha_{a,\varepsilon_0}}
{\partial y^\alpha_i}\zeta_i\Bigr) (\lambda^\alpha_1,\cdots, \lambda^\alpha_l)^t.$$
It follows that $\zeta=B(\lambda^\alpha_1,\cdots, \lambda^\alpha_l)^t$ for some
$B\in\mbox{\Bb R}$. If $\zeta\ne 0$ then it holds that
$$\frac{2A}{A-1}=t\sum^l_{i=1}\frac{\partial H^\alpha_{a,\varepsilon_0}}
{\partial y^\alpha_i}\lambda^\alpha_i.\eqno(2.18)$$
By (ii) the absolute value of the right hand of (2.18) is less than $\epsilon<1$,
and the left hand of it is more than 2 for every $A>1$. This contradiction
shows that $S$ is nonsingular.\hfill$\Box$\vspace{2mm}

Having this key lemma we can prove Proposition 1.4.B as follows.

\noindent{\it Proof of Proposition 1.4.B.}\hspace{3mm}
Without loss of generality we may assume $(Q,\sigma)$ to be exact.
By Lemma 2.4 there exists an
 open  neighbourhood ${\cal U}$ of $L$ in $Q$ and an embedding 
 $\varphi:{\cal U}\to T^\ast L$
such that $\varphi(L)$ coincides with the zero section of $T^\ast L$
and $\varphi^\ast(\omega^L_{\rm can}+ \pi^\ast_L(\sigma|_L))=\sigma$.
Applying Lemma 2.5 to $\lambda=\tau|_L$, $a=2c$ and $\epsilon=1$ we get a $b_1>0$
and a smooth function $K_{a,1}$. Let us take $\delta_0(c)>0$ so small that
the diffeomorphism $\Phi_\delta: T^\ast L\to T^\ast L$ given by
$$(m,v)\mapsto (m,\delta_0v)\eqno(2.19)$$
satisfies: 
$$\Phi_{\delta_0}(\{H_g\le 2b_1\})\subset \varphi({\cal U}).\eqno(2.20)$$
Notice that $\Phi_\delta$ is a symplectomorphis from $(\{H_g\le 2b_1\},
\omega^L_{\rm can}+ \pi^\ast_L(\sigma|_L))$ to
$$(\Phi_\delta(\{H_g\le 2b_1\}), \frac{1}{\delta}\omega^L_{\rm can}+ 
\pi^\ast_L(\sigma|_L))\eqno2.21)$$ 
for any $0<\delta\le\delta_0$. Denote by 
$\omega_{\delta, t}:=\frac{1}{\delta}\omega_{\rm can}^L +
\frac{1}{\delta} d(\pi_L^\ast\lambda) + 
t(1-\frac{1}{\delta})d(K_{a,1}\pi_L^\ast\lambda)$
for every $0<\delta\le\delta_0$ and $t\in [0,1]$.
By Lemma 2.5 (iii) they are symplectic forms on $T^\ast L$.
Moreover, $\omega_{\delta, 1}$ and $\frac{1}{\delta}\omega^L_{\rm can}+ 
\pi^\ast_L(\sigma|_L)$ are same on $\Phi_{\delta}(\{H_g\le 2c\})$.
By the construction of $\omega_{\delta, t}$ in lemma 2.5 we have
$$\varphi^\ast\omega_{\delta, t}=
\frac{1}{\delta}\sigma + t(1-\frac{1}{\delta})
d(\varphi^\ast(K^\alpha_{a,1}\pi^\ast_L(\tau|_L))).\eqno(2.22)$$
They are all symplectic forms on ${\cal U}$ and equal to 
$\frac{1}{\delta}\sigma$ near $\partial{\cal U}$.
Let us denote by
 $\widehat\omega_{\delta, t}=\frac{1}{\delta}\sigma + t(1-\frac{1}{\delta})
d(\varphi^\ast(K^\alpha_{a,1}\pi^\ast_L(\tau|_L)))$ for $t\in [0,1]$.
Since $\widehat\omega_{\delta, t}$ may naturally be extended onto symplectic
forms on $Q$ by assuming them being $\frac{1}{\delta}\sigma$
outside ${\cal U}$ we still denote them by $\widehat\omega_{\delta, t}$.
Using Moser's technique one can show that there exists a diffeomorphism
$F\in{\rm Diff}(Q)$ such that
$$F^\ast(\frac{1}{\delta}\sigma)=\widehat\omega_{\delta, 1}.\eqno(2.23)$$
Then the desired embedding $E_\delta$ is given by the composition
$$F\circ\varphi^{-1}\circ\Phi_\delta|_{\{H_g\le c\}}.\eqno(2.24)$$

The second claim follows from the first one and Lemma 2.2. 
For the third claim note that $F$ may be identity outside ${\cal U}$.
\hfill$\Box$\vspace{2mm}

\section{Proof of the main results}

\subsection{Proofs of Theorem 1.5.A and Corollary 1.5.B}

\noindent{\it Proof of Theorem 1.5.A}.\hspace{4mm}
Our idea of proof is to show that $C_{HZ}(\{ H_g\le c\}, \omega)$ is finite for
any $c>0$. For this purpose we first use Lemma 2.3 to get a symplectomorphism
$\Theta$ from $(T^\ast N,\omega)$ to 
$$(T^\ast N,\widehat\omega)=(T^\ast M, \omega^M_{\rm can}+\pi_M^\ast\widehat\Omega)
\times (T^\ast L, \omega^L_{\rm can}+\pi^\ast_L(\omega_0|_L)). \eqno(3.1)$$
Here $\widehat\omega$ and $\widehat\Omega$ are given by (2.7) (2.8) respectively.
Since the symplectic capacity is symplectic invariant we only need to prove that
$C_{HZ}(\{H_g\le c\},\widehat\omega)$ is finite for all $c>0$.
Fix a $c>0$ and take the Riemannian metric $g_1$ on $M$ and $g_2$ on $L$. Then
there exist positive numbers $c_1$ and $c_2$ such that
$\{z_1\in T^\ast M\,|\,H_{g_1}(z_1)\le c_1\}\times
\{z_2\in T^\ast L\,|\,H_{g_2}(z_2)\le c_2\}$ contains $\{H_g\le c\}$.
Therefore we only need to prove
$$C_{HZ}(\{H_{g_1}\le c_1\}\times\{H_{g_2}\le c_2\}, (\omega^M_{\rm can}+
\pi_M^\ast\widehat\Omega)\oplus (\omega^L_{\rm can}+\pi^\ast_L(\omega_0|_L)))<+
\infty.\eqno(3.2)$$
By Proposition 2.6 there exist a $\delta>0$ and a symplectic embedding of
codimension zero $E_\delta$ from  
$(\{z_2\in T^\ast L\,|\,H_{g_2}(z_2)\le c_2\}, \omega^L_{\rm can}+
\pi^\ast_L(\omega_0|_L))$ into $(\mbox{\Bb R}^{2l},\frac{1}{\delta}\omega_0)$.
We may assume the image of it to be contained in $Z^{2l}(R)$ for some
large $R>0$ since this image set is compact. Hence the monotonicity of
the symplectic capacity implies that the left hand of (3.2) is less than
\begin{eqnarray*}
&&C_{HZ}\bigl(\{z_1\in T^\ast M\,|\, H_{g_1}\le c_1\}\times Z^{2l}(R),
,(\omega^M_{\rm can}+\pi_M^\ast\widehat\Omega)\oplus \frac{1}{\delta}\omega_0\bigr)\\
&&< C_{HZ}\bigl(T^\ast M\times Z^{2l}(R), (\omega^M_{\rm can}+
\pi_M^\ast\widehat\Omega)\oplus \frac{1}{\delta}\omega_0\bigr)\\
&&\le\frac{\pi R^2}{\delta}<+\infty.
\end{eqnarray*}
Here in the final step we use Theorem 2.1 and fact that
$$m(T^\ast M, \omega^M_{\rm can}+\pi_M^\ast\widehat\Omega)=m(M,\widehat\Omega)=
m(M,\Omega)=+\infty.$$

If $L$ is simply connected it follows from the arguments above and remarks in [Lu2]
that $\bar C_{HZ}(T^\ast N,\omega)<+\infty$.
This completes proof of Theorem 1.5.A. \hfill$\Box$\vspace{2mm}

\noindent{\it Proof of Corollary 1.5.B}.\hspace{4mm}
Since $\pi_1(M)$ is a  finite group we choose a simply connected finite cover
$q_M:\widetilde M\to M$. Denote by $\widetilde N:=\widetilde M\times L$.
Notice that $H^2(L,\mbox{\Bb R})=0$. It directly follows from
the K\"unneth formula that $H^2(\widetilde N,\mbox{\Bb R})$ and
 $H^2(\widetilde M,\mbox{\Bb R})$ is isomorphic, and thus
$$\bigcup_{\phi\in{\rm Diff}(\widetilde N)}
\phi^\ast(P^\ast_{\widetilde M}(H^2({\widetilde M},\mbox{\Bb R})))=
H^2(\widetilde N,\mbox{\Bb R}).$$
That is, the requirement that the de Rham cohomology class of the related
closed $2$-form belongs to the set in (1.2) corresponding to 
$\widetilde N$ can always be satisfied.
 
Denote by 
${\widetilde q}_M:=q_M\times id_L$, 
 and 
$${\widetilde Q}_M:T^\ast\widetilde N\to T^\ast N:
(m,v)\mapsto ({\widetilde q}_M(m), [d{\widetilde q}_M(m)^{-1}]^\ast(v)).$$
For $\widetilde\Omega:={\widetilde q}_M^\ast\Omega$ it is easily checked that
$$\widetilde\omega:=\omega^{\widetilde N}_{\rm can}+
\pi_{\widetilde N}^\ast\widetilde\Omega=
{\widetilde Q}_M^\ast( \omega^N_{\rm can}+\pi_N^\ast\Omega)\quad{\rm and}
\quad m(\widetilde N,\widetilde\Omega)=m(N,\Omega).$$

Applying Theorem 1.5.A to the symplectic manifold
  $(T^\ast\widetilde N, \omega^{\widetilde N}_{\rm can}+
\pi_{\widetilde N}^\ast\widetilde\Omega)$ and the pullback metric 
$\tilde g:={\widetilde q}_M^\ast g$ we may get a  nonconstant Hamiltonian periodic
orbit on the generic levels $\{z\in T^\ast\widetilde N| H_{\tilde g}(z)=c\}$ of
 $X_{H_{\tilde g}}$ which is the Hamiltonian vector field of $H_{\widetilde g}$
with respect to the symplectic form $\widetilde\omega$.
Notice that the submersion  ${\widetilde Q}_M$ maps a nonconstant Hamiltonian 
periodic orbit of $X_{H_{\tilde g}}$    
 on $\{z\in T^\ast\widetilde N| H_{\tilde g}(z)=c\}$ to
a nonconstant one on $\{z\in T^\ast N| H_g(z)=c\}$ of $X_{H_g}$ with respect
to $\omega$. Corollary 1.5.B is proved. 
\hfill$\Box$\vspace{2mm}

\subsection{Proof of Theorem 1.5.C}

The ideas are similar to that of 
Theorem 1.5.A. Since $L$ is a closed Lagrangian submanifold
of $(\mbox{\Bb R}^{2l},\omega_0)$ one can directly use Weinstein's Lagrangian 
neighborhood theorem to get a symplectic embedding
 $\varphi:({\cal U},\omega_0)\to (T^\ast L,\omega^L_{\rm can})$
with $\varphi|_L=id$. Here ${\cal U}$ is an open neighborhood of $L$ in 
$\mbox{\Bb R}^{2l}$. As in the proof of Theorem 1.5.A we have a closed $2$-form
$\widehat\Omega$ on $M$ determined by  (2.8) and a symplectomorphism
$\Theta$ from $(T^\ast N,\omega)$ to
$$(T^\ast N,\widehat\omega):=(T^\ast M, \omega^M_{\rm can}+\pi_M^\ast\widehat\Omega)
\times (T^\ast L, \omega^L_{\rm can}).\eqno(3.3)$$
Moreover, we have also a symplectic embedding of codimension zero
$$\Upsilon:(T^\ast M,\omega^M_{\rm can}+\pi_M^\ast\widehat\Omega)\times
({\cal U},\omega_0)\to (T^\ast M,\omega^M_{\rm can}+\pi_M^\ast\widehat\Omega)\times
(T^\ast L, \omega^L_{\rm can})\eqno(3.4)$$
given by $(z_1, z_2)\mapsto (z_1,\varphi(z_2))$, whose the image  
is an open neighborhood of zero section of $T^\ast N$.  
For a given level $E_c=\{z\in T^\ast N\,|\, H_g(z)=c\}$ with $c>0$ 
 we can not guarantee that $\Theta(E_c)$ is contained in the image of $\Upsilon$ 
because $\Theta$ does not necessarily map the zero section to the zero section.
 But we can always take a  $\delta>0$ so small that the diffeomorphism
 $\Psi_{\delta}:T^\ast N\to T^\ast N$ given by
$$(m,v)=((m_1, m_2), (v_1, v_2))\mapsto ((m_1, m_2), (v_1, \delta v_2)),\eqno(3.5)$$
 maps $\Theta(E_c)$ into ${\rm Im}(\Upsilon)$. Let us denote by
$$\bar\delta=\bar\delta(c, g,\Theta,\Upsilon)>0\eqno(3.6)$$
the supreme of all  such $\delta>0$. 
Then for every $\delta\in(0,\bar\delta)$ it holds that
$$\Psi_{\delta}(\Theta(E_c))\subset{\rm Im}(\Upsilon).\eqno(3.7)$$
For such $\delta$, the composition
 $(\Upsilon|_{{\rm Im}(\Upsilon)})^{-1}\circ\Psi_{\delta}\circ\Theta$ is a 
 symplectic embedding of codimension zero from an open submanifold of
 $(T^\ast N,\omega)$ containing $E_c$ to 
 $$(T^\ast M, \omega^M_{\rm can}+\pi_M^\ast\widehat\Omega)\times 
({\cal U},\omega_0/\delta).\eqno(3.8)$$
 Denote by
 $$\Lambda(\Upsilon,\delta,\Theta,g,c):=
 P_{\cal U}[(\Upsilon|_{{\rm Im}(\Upsilon)})^{-1}\circ\Psi_{\delta}\circ
\Theta(E_c)],\eqno(3. 9)$$
 where $P_{\cal U}:T^\ast M\times{\cal U}\to {\cal U}$ is the natural projection.
It is a compact subset of open set ${\cal U}$ and is contained in
${\cal U}\cap Z^{2l}(r)$ for some $r>0$. Let us define 
$$r(\delta, c, g,\Theta,\Upsilon)\eqno(3.10)$$
the infimum of all $r>0$ such that
$$\Lambda(\Upsilon,\delta,\Theta,g,c)\subseteq{\cal U}\cap F(Z^{2l}(r))$$
for some $F\in{\rm Symp}(\mbox{\Bb R}^{2l},\omega_0)\}$. We also define
$$r(\bar\delta,c,g,\Theta,\Upsilon)=\inf_{0<\delta<\bar\delta}r(\delta, c,
 g,\Theta,\Upsilon)/\sqrt{\delta},\eqno(3.11)$$
where $Z^{2l}(r)$ is as in Theorem 2.1, then we easily prove that
$0< r(\bar\delta,c,g,\Theta,\Upsilon)<+\infty$ 
since ${\cal U}$ is a bounded open subset of $\mbox{\Bb R}^{2l}$.
Moreover, for each $R>r(\delta, c, g,\Theta,\Upsilon)$ there exists a 
$F\in{\rm Symp}(\mbox{\Bb R}^{2l},\omega_0)$ such that
$$\Lambda(\Upsilon,\delta,\Theta,g,c)\subseteq{\cal U}\cap F(Z^{2l}(R)).$$
 Furthermore, we define
$$r(c,g,\Omega):=\inf r(\bar\delta, c, g,\Theta, \Upsilon),\eqno(3.12)$$
where the infimum is taken over all possible $(\Theta, \Upsilon)$
satisfying the above arguments. Then the function
$$\Gamma^g_\Omega: [0, +\infty)\to [0,+\infty),\;c\mapsto r(c,g,\Omega)\eqno(3.13)$$
will satisfy the requirements of Theorem 1.5.C. To see these let
$c>0$ such that
$$\Gamma^g_\Omega(c)<\sqrt{m(N,\Omega)/\pi}.\eqno(3.14)$$
Then by (3.12) (3.13) we have
$$r(\bar\delta, c, g,\Theta, \Upsilon)<\sqrt{m(N,\Omega)/\pi}\eqno(3.15)$$
for some choice $(\Theta,\Upsilon)$,
and therefore from (3.11) it follows that there exists a $\delta\in (0,\bar\delta)$ 
such that
$$r(\delta, c, g,\Theta,\Upsilon)/\sqrt{\delta}<\sqrt{m(N,\Omega)/\pi}.
\eqno(3.16)$$
Let $\varepsilon>0$ satisfy
$$(r(\delta, c, g,\Theta,\Upsilon)+\varepsilon)/\sqrt{\delta}<
\sqrt{m(N,\Omega)/\pi}.
\eqno(3.17)$$
Then by the definition of $r(\delta, c, g,\Theta,\Upsilon)$ in (3.10)
there exists a $F\in{\rm Symp}(\mbox{\Bb R}^{2l},\omega_0)$ such that
$$\Lambda(\Upsilon,\delta,\Theta,g,c)
\subseteq{\cal U}\cap F(Z^{2l}(r(\delta, c, g,\Theta,\Upsilon)+
\varepsilon)).$$
Note that the left side is a compact subset and the right side is an open
set.  This implies that for sufficiently small $\epsilon>0$ 
$$P_{\cal U}[(\Upsilon|_{{\rm Im}(\Upsilon)})^{-1}\circ\Psi_{\delta}\circ
\Theta(U(g,c,\epsilon))]\subseteq{\cal U}\cap F(Z^{2l}(r(\delta, c, g,
\Theta,\Upsilon)+\varepsilon)).$$
Hence $\Upsilon|_{{\rm Im}(\Upsilon)})^{-1}\circ\Psi_{\delta}\circ\Theta
(U(g,c,\epsilon))$ is contained in
$T^\ast M\times {\cal U}\cap F(Z^{2l}(r(\delta, c, g,\Theta,\Upsilon)+\varepsilon))$.
This implies that the map
$(id\times F)^{-1}\circ(\Upsilon|_{{\rm Im}(\Upsilon)})^{-1}
\circ\Psi_{\delta}\circ\Theta$
 symplectically embeds $(U(g,c,\epsilon), \omega)$ into
$$(T^\ast M, \omega^M_{\rm can}+\pi_M^\ast\widehat\Omega)\times 
(Z^{2l}(r(\delta, c, g,\Theta,\Upsilon)+\varepsilon),
\omega_0/\delta).\eqno(3.18)$$
Using Theorem 2.1 and the fact that  (3.18) is symplectomorphic to
$$(T^\ast M, \omega^M_{\rm can}+\pi_M^\ast\widehat\Omega)\times 
(Z^{2l}([r(\delta, c, g,\Theta,\Upsilon)+\varepsilon]/\sqrt{\delta}),
\omega_0),\eqno(3.19)$$
we obtain that
$$C_{HZ}(U(g,c,\epsilon), \omega)\le\pi \Bigl([r(\delta, c, g,
\Theta,\Upsilon)+\varepsilon]/\sqrt{\delta}\Bigr)^2.$$
Hence (3.17) gives
$$\liminf_{\epsilon\to 0+}C_{HZ}(U(g,c,\epsilon), \omega)< m(N,\Omega).$$
The monotonicity of symplectic capacity $C_{HZ}$ leads to (1.3) directly.

 Finally, we prove that the function $\Gamma^g_\Omega$ is upper semi-continuous.
To this goal we only need to prove that the function 
$c\mapsto r(\delta, c, g,\Theta,\Upsilon)$
defined in (3.19) is upper semi-continuous.
 Fix a $c>0$ and a real number
$\lambda > r(\delta, c, g,\Theta,\Upsilon)$, we wish to prove
that if $c^\prime>0$ is sufficiently close to $c$ then
$$\lambda > r(\delta, c^\prime, g,\Theta,\Upsilon).$$
Otherwise, assume that there exists a sequence of $c_n>0$ such that
$$c_n\to c\;(n\to\infty)\qquad{\rm and}\qquad
\lambda \le r(\delta, c_n, g,\Theta,\Upsilon).\eqno(3.20)$$
Taking sufficiently small $\epsilon>0$ such that
$$\lambda-2\epsilon > r(\delta, c, g,\Theta,\Upsilon),$$
then by (3.10) we have
$$\Lambda(\Upsilon,\delta, \Theta, g, c_n)\not\subseteq{\cal U}\cap F(Z^{2l}(
\lambda-\epsilon)),\;\;\forall F\in{\rm Symp}(\mbox{\Bb R}^{2l},\omega_0).$$
On the other hand when $n\to\infty$ the compact subsets
$\Lambda(\Upsilon,\delta,\Theta, g, c_n)$ converges to the compact subset
$\Lambda(\Upsilon,\delta,\Theta, g, c)$ in the Hausdorff metric (even stronger sense).
Hence for every open neighborhood ${\cal V}$ of $\Lambda(\Upsilon,\delta,\Theta, g,c)$
the sets $\Lambda(\Upsilon,\delta,\Theta, g,c_n)$ can be contained in ${\cal V}$ for
sufficiently large $n$. If we understand the cyclinder $Z^{2n}(r)$ in Theorem 2.1 
 as the open cyclinder, then for any fixed
$F\in{\rm Symp}(\mbox{\Bb R}^{2l},\omega_0)$ the set
${\cal U}\cap F(Z^{2l}(\lambda-2\epsilon))$ is an open neighborhood of
$\Lambda(\Upsilon,\delta,\Theta, g, c)$ and thus
$$\Lambda(\Upsilon,\delta,\Theta, g,c_n)\subseteq{\cal U}\cap
 F(Z^{2l}(\lambda-2\epsilon))$$
for sufficiently large $n$. This shows that
 $$r(\delta, c_n, g,\Theta,\Upsilon)\le\lambda-2\epsilon,$$
which contradicts (3.20).

By the problem F.(d) on the page 101 of [K] it is easy to know that
the function $c\mapsto r(\bar\delta,c, g,\Theta,\Upsilon)$ and thus
$c\mapsto r(c, g, \Omega)$ are upper semi-continous.
The proof of Theorem 1.5.C is completed.
\hfill$\Box$\vspace{2mm}

\noindent{\bf Remark 3.2.A}\hspace{2mm}From the proof of Theorem 1.5.C we may see
 that the condition $\Omega|_{\pi_2(N)}=0$ may be weakened to the case that 
$m(N,\Omega)$ is only finite positive number. We here do not pursue it.

\subsection{Proofs of Theorem 1.5.E and Corollary 1.5.F}

\noindent{\it Proof of Theorem 1.5.E}\hspace{2mm}
Let $c>0$ such that $\Xi^g_\Omega(c)< m(N,\Omega)$.
For any $\varepsilon>0$ satisfying $\Xi^g_\Omega(c) + 2\varepsilon< m(N,\Omega)$,
by (1.5), we have $(\Theta,\widehat\Omega)$ and $\epsilon>0$
such that
$$\Bigl|\int_{\bar P_S(\Theta((U(g,c,\epsilon)))}
\omega^S_{\rm can}\Bigr|<\Xi^g_\Omega(c)+\varepsilon.\eqno(3.21)$$
Since $\bar P_S(\Theta((U(g,c,\epsilon)))\subset\mbox{\Bb S}^1\times\mbox{\Bb R}$
and the symplectomorphisms on $2$-dimensional symplectic manifolds are equivalient to
the diffeomorphisms presvering area we may find a symplectic embedding $F$ from
$\bar P_S(\Theta((U(g,c,\epsilon)))$ into a disk $B^2$ of area 
 $\Xi^g_\Omega(c)+ 2\varepsilon$ centred at origin in
$\mbox{\Bb R}^2$.  Then $id\times F$ symplectically embeds $(U(g,c,\epsilon),\omega)$
into $(T^\ast M,\omega^M_{\rm can}+\pi^\ast_M\widehat\Omega)\times (B^2,\omega_0)$.
Thus  by Theorem 2.1 it holds that
$$C_{HZ}(U(g,c,\epsilon),\omega)\le \Xi^g_\Omega(c)+ 2\varepsilon< m(N,\Omega).$$
Using the same reason as in the proof of Theorem 1.5.C one can get conclusions.

In order to prove (1.6), we denote by
$\Theta^{\alpha}_{\widehat\Omega}$ and $\Theta^{\lambda\alpha}_{\lambda\widehat\Omega}$
the corresponding symplectomorphisms to $(\Omega,\widehat\Omega,\alpha)$ and
$(\lambda\Omega,\lambda\widehat\Omega,\lambda\alpha)$ constructed in (2.8) and (2.12)
respectively, then it is easily checked that
\begin{eqnarray*}
&&\Theta^\alpha_\Omega(U(g, c,\epsilon))=\bigl\{\bigl(\phi^{-1}(m), d\phi(m)^\ast (v)
+\alpha(\phi^{-1}(m))\bigr)\,\bigm|\, (m, v)\in U(g, c,\epsilon)\bigr\},\\
&&\Theta^{\lambda\alpha}_{\lambda\Omega}(U(g,\lambda c,\lambda\epsilon))=\bigl\{
\bigl(\phi^{-1}(m), \lambda [d\phi(m)^\ast (v)+
\alpha(\phi^{-1}(m))]\bigr)\,\bigm|\, (m, v)\in U(g, c,\epsilon)\bigr\}
\end{eqnarray*}
and thus
$$\int_{\bar P_S(\Theta^{\lambda\alpha}_{\lambda\Omega}(U(g, \lambda c,
\lambda\epsilon)))}\omega^S_{\rm can}
=\lambda\int_{\bar P_S(\Theta(U(g,c,\epsilon)))}\omega^S_{\rm can}.$$
This can lead to (1.6). 
The upper semi-continousity of the function $\Xi^H_\Omega$ may be proved similarily
 as in Theorem 1.5.C.
\hfill$\Box$\vspace{2mm}

\noindent{\it Proof of Corollary 1.5.F}\hspace{2mm}
Under the assumptions of Corollary 1.5.F we have a closed two-form $\widehat\Omega_1$
on $M$ with $m(N,\Omega)=m(M,\widehat\Omega_1)$ and a diffeomorphism
$\Theta_1\in{\rm Diff}(T^\ast N)$ such that
$\omega=\Theta_1^\ast\widehat\omega_1$, where
$$\widehat\omega_1:=\omega_{\rm can} + \pi_N^\ast(P_M^\ast\widehat\Omega_1).$$
Writting $\overline M:=M\times T^{n-1}$ and $N=\overline M\times\mbox{\Bb S}^1$,
it is easily checked that the closed two-form $\Omega_2:=P_M^\ast\widehat\Omega_1$ 
on $N$ satisfies the requirements of Theorem 1.5.E. For 
$\widehat\omega_1=\omega_{\rm can} + \pi_N^\ast{\widehat\Omega}_2$ we can obtain
a closed two-form ${\widehat\Omega}_3$ on $\overline M$ and a  diffeomorphism
$\Theta_2\in{\rm Diff}(T^\ast N)$ such that
$\Theta_2{\widehat\omega}_2^\ast={\widehat\omega}_1$, where
$${\widehat\omega}_2:=\omega_{\rm can} + 
\pi_N^\ast(P_{\overline M}^\ast{\widehat\Omega}_2).$$
Denote by $\Theta=\Theta_1\circ\Theta_2$ and set
$${\widehat\Xi}^g_{\Omega}(c):={\inf}_{(\Theta,{\widehat\Omega}_1,
{\widehat\Omega}_2, \epsilon)}
\Bigl|\int_{\bar P_S(\Theta((U(g,c,\epsilon)))}
\omega^S_{\rm can}\Bigr|,\eqno(3.22)$$
where the infimum is taken over all triples 
$(\Theta, {\widehat\Omega}_1,{\widehat\Omega}_2)$
satisfying the above arguments and all $\epsilon>0$.  
Then it is easily proved that for every $c>0$ 
$$\lim_{\epsilon\to 0+}C_{HZ}(U(g,c,\epsilon),\omega)\le {\widehat\Xi}^g_\Omega(c)< 
m(N,\Omega)\eqno(3.23)$$
if ${\widehat\Xi}^g_\Omega(c)< m(N,\Omega)$, where we use that 
$m(N,\Omega)=m(M,{\widehat\Omega}_1)=m(\overline M,{\widehat\Omega}_2)$.
Specially, similar to the proof of (1.6) in Theorem 1.5.E 
${\widehat\Xi}^{g}_\Omega$ also satisfies:
$${\widehat\Xi}^g_{\lambda\Omega}(\lambda c)=\lambda {\widehat\Xi}^g_\Omega(c). 
\eqno(3.24)$$
for all $\lambda>0$ and $c>0$. The proof is completed.
\hfill$\Box$\vspace{2mm} 

\subsection{Proof of Theorem 1.5.I}

\noindent{\bf Case 1}.\hspace{2mm}$n>2$.
 
Notice that 
$T^\ast T^n=T^n\times\mbox{\Bb R}^n$. We may denote  by $(x_1,\cdots,x_n; y_1,
\cdots, y_n)$ the coordinate in it. For the sake of clearness we give some reductions,
which are, either more or less, contained in [Gi1][GiK][J2]. 
Since $[dx_i\wedge dx_j]$( $1\le i<j\le n$) form a basis of vector space
 $H_{de}^2(T^n;\mbox{\Bb R})$, and every
closed two-form $\Omega$ on $T^n$ must have the following form:
$$\Omega=\sum_{i,j}q_{ij}(x)dx_i\wedge dx_j,$$
where the smooth function $q_{ij}(x)=-q_{ji}(x)$ are $1$-periodic for
each variable $x_i$. From them it follows that there exist
constants $c_{ij}$ such that
$$[\Omega]=\sum_{i<j}c_{ij}[dx_i\wedge dx_j].$$  
In fact $c_{ij}=\int_{T^n}q_{ij}(x)$.
Setting $b_{ij}=\frac{1}{2}c_{ij}$ then there exists a $1$-form 
$\alpha$ on $T^n$ such that
$$\Omega=\sum_{i,j}b_{ij}dx_i\wedge dx_j + d\alpha.$$
Moreover, $[\Omega]\ne 0$ if and only if there at least exists a $b_{ij}\ne 0$.
 We may write $\alpha$ as 
 $$\alpha=\sum^n_{i=1}a_i(x)dx_i,$$ 
 where the smooth functions $a_i(x)$ are $1$-periodic for each variable
 $x_i$. It is easy to see that  the transformation
 $\psi:T^n\times\mbox{\Bb R}^n\to T^n\times\mbox{\Bb R}^n:
 (x,y)\mapsto (X,Y)=(x, y-a(x))$ is a symplectmorphism from
$(T^n\times\mbox{\Bb R}^n, \omega)$ to
$(T^n\times\mbox{\Bb R}^n, \sum_{i=1}^n dX_i\wedge dY_i +
\sum_{i,j}b_{ij}dX_i\wedge dX_j)$, where
 $\omega:=\sum_{i=1}^n dx_i\wedge dy_i +\Omega$.
Now  there at least exists  a  $b_{ij}\ne 0$. It is this condition which leads to
one to be able  prove that $(T^n\times\mbox{\Bb R}^n  ,\omega)$ is symplectomorphic 
to the product  $(\mbox{\Bb R}^{2k}\times W_1, \omega_0\oplus\sigma)$ with
$k\ge 1$([GiK]), where $\omega_0$ is the standard symplectic form on 
$\mbox{\Bb R}^{2k}$  and $\sigma$ is a  translation-invariant symplectic form on 
$W_1=\mbox{\Bb R}^{n-2k}\times T^{n}$ . Note that here the assumption $n>2$
is used. 
Let us denote the symplectomorphism by
$\phi$. For a given metric $g$ on $T^n$ and $c>0$ we also denote by
$$U_c=\{(x,y)\in T^n\times\mbox{\Bb R}^n\,|\, g(x)(y, y)\le c^2\},$$  
then there exists a $r_c>0$ such that
$$\phi(U_c)\subset W_1\times B^{2r}(r_c).$$
Denote by $V_c= W_1\times B^{2r}(r_c)$ and the inclusion maps
\begin{eqnarray*}
&I_{U_c}: U_c\hookrightarrow T^n\times\mbox{\Bb R}^n,\quad &I_{\phi(U_c)}:
\phi(U_c)\hookrightarrow\mbox{\Bb R}^{2k}\times W_1\\
&I_{\phi(U_c)V_c}: \phi(U_c)\hookrightarrow V_c,\quad &I_{V_c}:
V_c\hookrightarrow\mbox{\Bb R}^{2k}\times W_1.
\end{eqnarray*}
We have the following commutative diagram:

\begin{center}\setlength{\unitlength}{1mm}
\begin{picture}(80,30)
\thinlines
\put(8,25){$T^n\times\mbox{\Bb R}^n$}
\put(25,25){\vector(1,0){44}}
\put(45,27){$\phi$}
\put(15,22){\vector(0,-1){16}}
\put(8,13){$I_{U_c}$}
\put(70,25){$\mbox{\Bb R}^{2k}\times W_1$}
\put(80,22){\vector(0,-1){16}}
\put(70,13){$I_{\phi(U_c)}$}
\put(13,0){$U_c$}
\put(23,1){\vector(1,0){47}}
\put(43,3){$\phi|_{U_c}$}
\put(75,0){$\phi(U_c)$}
\end{picture}\end{center}

Then the induced homomorphisms among their first homotopy groups satisfy:
$$\phi_\ast\circ I_{U_c\ast}=I_{\phi(U_c)\ast}\circ (\phi|_{U_c})_\ast.$$  
Since $\phi_\ast$,  $I_{U_c\ast}$ and  $(\phi|_{U_c})_\ast$ are all isomorphisms
the homomorphism $I_{\phi(U_c)\ast}$ is also an isomorphism. But
$I_{\phi(U_c)}=I_{V_c}\circ I_{\phi(U_c)V_c}$. We get that
$I_{\phi(U_c)\ast}=I_{V_c\ast}\circ I_{\phi(U_c)V_c\ast}$.
This implies that $I_{\phi(U_c)V_c\ast}: \pi_1(\phi(U_c))\to 
\pi_1(V_c)$ must be injective.
 Using Theorem 2.1 and the weak monotonicity of $\bar C_{HZ}$ we have:
$$\bar C_{HZ}(U_c,\omega)\le \bar C_{HZ}(V_c,\omega_0\oplus\sigma)\le
\bar C_{HZ}(B^{2k}(r_c)\times W_1, \omega_0\oplus\sigma)\le \pi r_c^2.$$
Hence, for generic $c>0$ the level $E_c$ carries a nonconstant periodic
orbit $z=z(t)$ of $X_{H_g}$, which is contractible in $U_c$. 
Since the fibre projection from $U_c$ to $T^n$ induces an isomorphism 
$\pi_1(U_c)\to\pi_1(T^n)$ the projection of $z=z(t)$ to the base $T^n$
is contractible. This leads to our claim.

\noindent{\bf Case 2}.\hspace{2mm}$n=2$.

Denote by $\sigma_0$ the standard symplectic form on $T^2$.  Since $[\Omega]\ne 0$
there exists a constant $c_0\ne 0$ such that $[\Omega]=c_0[\sigma_0]$.
In particular, Moser theorem shows that $(T^2,\Omega)$ is symplectomorphic to
$(T^2, c_0\sigma_0)$.  Notice that there exists a symplectomorphism
from $(T^\ast T^2,\omega)$ to $(T^2\times\mbox{\Bb R}^2, \Omega\oplus\omega_0)$
which maps the zero section to $T^2\times\{0\}$ for the standard symplectic form 
$\omega_0$ on $\mbox{\Bb R}^2$. They together must induce
such a symplectomorphism $\Psi$ from $(T^\ast T^2,\omega)$ to 
$(T^2\times\mbox{\Bb R}^2, (c_0\sigma_0)\oplus\omega_0)$.
Now for any $U_c$ as above there exists a $r_c>0$ such that
$\Psi(U_c)\subset T^2\times B^2(r_c)$. It is easily checked that
$\Psi$ also induces an injective map 
$\Psi_\ast:\pi_1(U_c)\to\pi_1(T^2\times B^2(r_c))$. Hence 
$\bar C_{HZ}(U_c,\omega)\le\bar C_{HZ}(T^2\times B^2(r_c), 
(c_0\sigma_0)\oplus\omega_0)\le\pi r_c^2$. This leads to the conclusion again.
\hfill$\Box$\vspace{2mm} 

\section{ The concluding remarks}

  Our methods can actually be used to deal with a more general question than PMMQ 
above.

\noindent{\bf Definition 4.1.}\quad {\it A smooth, bounded from below functiom
$H:M\to\mbox{\Bb R}$ is called  {\bf strong proper} if
the sublevel $\{z\in M\,|\, H(z)\le c\}$ is compact for every $c\in{\rm Im}(H)$.}

Clearly, $H_g$ is a strong proper on $T^\ast N$. Let $H:T^\ast N\to\mbox{\Bb R}$ 
a strong proper function and
$$\omega:=\omega_{\rm can}+\Omega\eqno(4.1)$$
is a symplectic form  on $T^\ast N$, where $\Omega$ is a closed $2$-form on
$T^\ast N$. One may ask the following more general question:

\noindent{\bf Question 4.2.}\quad Does the Hamiltonian flow of $H$ with respect
to the symplectic form in (4.1) has a nonconstant periodic orbit on the level
$\{H=c\}$  for every $c\in{\rm Im}(H)$?

The proof of Theorem 1.5.A can suitably be modified to show
 
 \noindent{\bf Theorem 4.3}\hspace{2mm}{\it Let $N=M\times L$ be in Theorem 1.5.A. 
Suppose that a closed two-form $\Omega$ on $T^\ast N$  such that
 \begin{description}
 \item[(i)] the $2$-form $\omega:=\omega_{\rm can}+\Omega$  is a symplectic form on
            $T^\ast N$ and $\Omega|_{\pi_2(T^\ast N)}=0$;
 \item[(ii)] there is a ${\cal A}\in{\rm Symp}(T^\ast N,\omega_{\rm can})$ such that
$[({\cal A}^\ast\Omega)|_N]$  belongs to the set in (1.2), and
 $(T^\ast N, \omega_{\rm can}+ {\cal A}^\ast\Omega)$ is
 symplectomorphic to $(T^\ast N, \omega_{\rm can}+
\pi_N^\ast(({\cal A}^\ast\Omega)|_N))$.
  \end{description}
 Then for every  strong proper smooth function $H:T^\ast N\to\mbox{\Bb R}$ 
 and every $c\in{\rm Im}(H)$ it holds that
 $$C_{HZ}(\{H\le c\},\omega)<+\infty.$$
 Consequently, for almost all $c\in {\rm Im}(H)$ the levels $\{H=c\}$ carries a 
nonconstant periodic  orbit of $X_H$, where $X_H$ is contraction by
 $i_{X_H}\omega=dH$. Furthermore, if $L$ is simply connected, and an regular value
 $c_0\in{\rm Im}(H)$ is such that the inclusion  
$\{H\le c_0\}\hookrightarrow T^\ast N$ induces an injective homomorphism
$\pi_1(\{H\le c_0\})\to\pi_1(T^\ast N)$ then for generic $c$ near $c_0$
 the levels $\{H=c\}$ carries a nonconstant periodic  orbit of $X_H$ with the
 contractible projection to $N$.}
 
For other theorems and corollaries in this paper similar results may also be
obtained. We omit them.

\noindent{\bf Remark 4.4.}\hspace{2mm}One may feel that checking the condition (ii) 
in Theorem 4.3 is difficult. But under some cases the first claim of it implies 
the second one. In fact, since $\pi_N:T^\ast N\to N$ induces an 
isomorphism $\pi_N^\ast:H^2(N,\mbox{\Bb R})\to H^2(T^\ast N,\mbox{\Bb R})$ 
we have  $[{\cal A}^\ast\Omega]=[\pi_N^\ast(({\cal A}^\ast\Omega)|_N)]$, and thus
$${\cal A}^\ast\Omega-\pi_N^\ast(({\cal A}^\ast\Omega)|_N)=d\alpha$$
for some one-form $\alpha$ on $T^\ast N$. Using Moser's technique one can prove that
 if all two-forms 
 $$\omega_{\rm can}+\pi_N^\ast(({\cal A}^\ast\Omega)|_N)+ t d\alpha,\;0\le t\le 1,$$
 are  symplectic forms, and $\alpha$ satisfies some conditions( for example,
 the norm $|\alpha|$ with respect to some complete metric on $T^\ast N$ is bounded), 
  then the second claim in (ii) may be satisfied.

\end{document}